\documentclass[10pt]{amsart}
\usepackage{amsmath,amssymb,amsthm}
\usepackage[all]{xy}
\usepackage{hyperref}

\theoremstyle{plain}
\newtheorem{The}{Theorem}[section]
\newtheorem{prop}[The]{Proposition}

\newtheorem{lem}[The]{Lemma}

\theoremstyle{definition}

\newtheorem{Ex}[The]{Example}
\newtheorem{Rem}[The]{Remark}

\newtheorem{Nots}[The]{Notations}

\newtheorem{Hyp}[The]{Assumption}

\newcommand{\C}{\mathbf{C}}

\newcommand{\A}{\mathbf{A}}

\newcommand{\PP}{\mathbf{P}}
\newcommand{\ui}{\textbf{i}}

\newcommand{\pr}{\textit{Proof.}\ }

\newcommand{\Ot}{\mathcal{O}_{\tau}}
\newcommand{\Ut}{U_{\tau}^-}
\newcommand{\U}{\mathcal{U}}
\renewcommand{\t}{\tau}
\renewcommand{\o}{\omega}
\renewcommand{\a}{\alpha}
\renewcommand{\b}{\beta}
\newcommand{\g}{\gamma}

\newcommand{\mult}{\operatorname{mult}}
\newcommand{\Sing}{\operatorname{Sing}}

\newcommand{\im}{\operatorname{Im}}

\renewcommand{\phi}{\varphi}

\begin{document}

\title[Multiplicity on a Richardson variety]{Multiplicity on a Richardson variety in a cominuscule $G/P$}
\author{Micha\"el Balan}
\date{April 21, 2011}
\address{Universit\'e de Valenciennes\\Laboratoire de Math\'ematiques\\Le Mont Houy -- ISTV2\\F-59313 Valenciennes Cedex 9\\France}
\email{michael.balan@univ-valenciennes.fr}

\begin{abstract}
We show that in a cominuscule partial flag variety $G/P$, the multiplicity of an arbitrary point on a Richardson
variety $X_{w}^{v}=X_w \cap X^v \subset G/P$ is the product of its multiplicities on the Schubert varieties
$X_w$ and $X^v$.
\end{abstract}

\maketitle

\section*{Introduction}

Richardson varieties, named after \cite{R}, are intersections of a Schubert variety and an opposite Schubert
variety inside a partial flag variety $G/P$ ($G$ a connected complex semi-simple group, $P$ a parabolic subgroup). They 
previously appeared in \cite[Ch.~XIV, \S4]{HP} and \cite{St}, as well as the corresponding open cells in \cite{De}.
They have since played a role in different contexts, such as equivariant K-theory \cite{LL}, positivity in Grothendieck
groups \cite{B2}, standard monomial theory \cite{BL}, Poisson geometry \cite{GY}, positroid varieties \cite{KLS}, 
and their generalizations \cite{KLS2,BC}.

On the other hand, singularities of Schubert varieties have been extensively studied in the last decades.
The singular locus of Schubert varieties in Grassmannians has been determined independently in \cite{Sv}
and \cite{La}, and more generally in a minuscule $G/P$ in \cite{LW}. In the
full flag variety of type~A$_n$, it has been determined independently in \cite{BW}, \cite{Co},
\cite{KLR}, and \cite{Ma}.

Moreover, the multiplicity of a singular point on a Schubert variety is known in several cases: when $G/P$ is minuscule of arbitrary type, 
or cominuscule of type~C$_n$, a recursive formula was given in \cite{LW}. A direct determinantal formula was given in \cite{RZ} for $G/P$ a
Grassmannian; it has been subsequently interpreted
in terms of non-intersecting lattice paths \cite{Kr}. The multiplicity problem has also been studied in relationship with Hilbert functions and 
Gr\"obner degenerations \cite{GR,KR,Kr2,KL2,RU2,RU}, as well as 
with $T$--equivariant cohomology \cite{I,IN,KT,K1,K2,LRS}. The problem of determining the multiplicity of a point in a Schubert variety in the 
full flag variety is more complicated; see \cite{W,LY,WY,WY2}. 

For Richardson varieties in a minuscule $G/P$, the multiplicity of a $T$-fixed point 
($T \subset P$ a maximal torus in $G$) has been determined by Kreiman and Lakshmibai \cite{KL} (for the Gr\"obner point of view, see also \cite{K} in
type A$_n$ and \cite{U} in orthogonal types).

In this paper, we determine the multiplicity of an arbitrary point\footnote{Note
that unlike in the case of a Schubert variety, this no longer follows from information about $T$-fixed points, as pointed
out in the introductions of \cite{KL} and \cite{K}.} on a Richardson variety in a cominuscule $G/P$.

Before stating the main result, let us fix some notation. Let $G,P,T$ be as above, with $G$ adjoint.
Let $X(T)$ be the character group of $T$, $R \subset X(T)$ the root system, and $W=N_G(T)/T$ its Weyl group.
Let $B\subset G$ be a Borel subgroup such that $T\subset B\subset P$: it determines a system of positive roots
$R^+$ and a system of simple roots $S$. Denote by $B^-$ the opposite Borel subgroup (\textit{i.e.}\ such that
$B\cap B^-=T$).

Let $W_P\subset W$ be the subgroup associated to $P$ (so that $W_G=W$ and $W_B$ is the trivial subgroup).
In the quotient $W^P=W/W_P$, every coset $wP$ contains a unique minimal element
for the Bruhat order $\leq$ on $W$, so we shall identify $W^P$ with the set of minimal representatives. The
$B$-orbit (resp.\ the $B^-$-orbit) of a $T$-fixed point $e_{\tau}=\tau P$ is called a Schubert cell (resp.\ an
opposite Schubert cell) in $G/P$, and denoted by $C_{\tau}$ (resp.\ $C^{\tau}$). Its closure is the Schubert variety
$X_\tau$ (resp.\ the opposite Schubert variety $X^{\tau}$).

If $v,w\in W^P$, then the intersection $X_w^v=X_w \cap X^v$ is called a \textit{Richardson variety}; it is non-empty
if and only if $v\le w$ (note that Schubert varieties are the particular cases $X_w=X_w^e$ and $X^v=X_{w_0}^v$, where $e,
w_0 \in W$ are the identity and the longest element, respectively).

Now assume $P$ to be maximal, and let $\a$ be the associated simple root (so that $W_P$ is generated by the reflections $s_\delta$ with $\delta\in S\setminus\{\alpha\}$).
Then $P$ (or $\a$) is said to be
\begin{itemize}
\item cominuscule if $\alpha$ occurs with a coefficient~1 in the decomposition
of the highest root of $R^+$;
 \item minuscule if $\a^\vee$ is cominuscule in the dual root system $R^\vee$.
\end{itemize}~\\

The main result of this paper is the following

\begin{The}\label{theformula} Assume $P$ is cominuscule. Let $m \in X_w^v$ be arbitrary, and denote by $\mu_w$ (resp.\ $\mu^v$, $\mu_w^v$) the
multiplicity of $m$ on $X_w$ (resp.\ $X^v$, $X_w^v$). Then
\begin{equation}\label{multmult}
\mu_w^v=\mu_w\ \mu^v.
\end{equation}
\end{The}

This result indeed determines the multiplicities on $X_w^v$, since those on $X_w$ and $X^v$ are known:
types A$_n$, D$_n$, E$_6$, E$_7$ are covered by \cite{LW}, Section~3 (since cominuscule is equivalent to minuscule in those
types), and type C$_n$ is covered by \cite{LW}, Section~4. The only remaining case, in type B$_n$ (\textit{cf.}\ the table below),
is elementary, and covered in the Appendix of the present paper for the sake of completeness.

Note that (\ref{multmult}) is exactly the result obtained in \cite{KL} for a $T$-fixed point in a minuscule $G/P$.\\

To prove the theorem, we shall use a description of the multiplicity using a central projection: namely, given a
projective variety $X\subset \PP^N$ and a point $m \in X$, we consider the projection $p_m$, of centre $m$, onto
a hyperplane not containing $m$. Then the multiplicity of $m$ on $X$ is the difference between the degree of $X$
and the projective degree of $p_m$. Note that the projective degree of $p_m$ is zero when $X$ is a cone. We
apply this description for $X$ the projective closure of the affine trace $X_w^v \cap \mathcal{O}_{\tau}$, where
$\Ot$ is an affine open subset of $G/P$ identified with $\A^N$. One then needs to know whether
the affine traces of $X_w,X_v,X_w^v$ are cones or not. In this setting, we can explain why we assume that $P$ is cominuscule:
\begin{itemize}
\item it implies that $X_w\cap\Ot$ is a cone over \emph{any} point of the cell $C_\tau$ (though this may not be the case for $X^v\cap\Ot$);
\item we relate the central projection $p_m$ to a map which turns out to be a $\C$-action if $P$ is cominuscule. It
is this $\C$-action which allows to prove all the necessary properties for $p_m$.
\end{itemize}

In Section~\ref{coordinates}, we give a system of local coordinates in which $X_w\cap\Ot$ is a cone over both
$e_\tau$ and $m$, and $X^v\cap\Ot$ over $e_\tau$. In Section~\ref{projection}, we prove Theorem~\ref{theformula}
assuming certain formulas for the degrees involved and that $X^v\cap\Ot$ is not a cone over $m$. These
assumptions are summarized in Proposition~\ref{theproposition}, and proved in Sections~\ref{proof}~and \ref{secinter}. 
The proofs are based on a $\C$-action linking the central projections of centres $m$ and
$e_\tau$; this action is defined and studied in Section~\ref{Caction}.

For the convenience of the reader, we give the minuscule and cominuscule weights in the following table:\\

\begin{tabular}{|c|c|}
\hline $A_n$ &
\begin{picture}(20,2)(-0.2,0)
\put(0,1){\circle*{0.4}} \put(0,1){\circle{0.6}} \put(0.2,1){\line(1,0){3}}
\put(3.4,1){\circle*{0.4}}\put(3.4,1){\circle{0.6}}\put(3.6,1){\line(1,0){3}}\put(6.8,1){\circle*{0.4}}
\put(6.8,1){\circle{0.6}} \put(7,1){\line(1,0){1}} \put(9,1){\circle*{0.2}} \put(10,1){\circle*{0.2}}
\put(11,1){\circle*{0.2}} \put(12,1){\line(1,0){1}} \put(13.2,1){\circle*{0.4}} \put(13.2,1){\circle{0.6}}
\put(13.4,1){\line(1,0){3}} \put(16.6,1){\circle*{0.4}} \put(16.6,1){\circle{0.6}} \put(0,0){1} \put(3.4,0){2}
\put(6.8,0){3} \put(13.2,0){$n-1$} \put(16.6,0){$n$}
\end{picture}\\
\hline $B_n$ &
\begin{picture}(20,2)(-0.2,0)
\put(0,1){\circle{0.4}} \put(0,1){\circle{0.6}} \put(0.2,1){\line(1,0){3}}
\put(3.4,1){\circle{0.4}}\put(3.6,1){\line(1,0){3}}\put(6.8,1){\circle{0.4}}\put(7,1){\line(1,0){1}}
\put(9,1){\circle*{0.2}} \put(10,1){\circle*{0.2}} \put(11,1){\circle*{0.2}} \put(12,1){\line(1,0){1}}
\put(13.2,1){\circle{0.4}}\put(13.4,1.1){\line(1,0){3}}\put(13.4,0.9){\line(1,0){3}}\put(15.6,1){\line(-2,1){1}}
\put(15.6,1){\line(-2,-1){1}}\put(16.6,1){\circle*{0.4}} \put(0,0){1} \put(3.4,0){2} \put(6.8,0){3}
\put(13.2,0){$n-1$} \put(16.6,0){$n$}
\end{picture}\\
\hline $C_n$ &
\begin{picture}(20,2)(-0.2,0)
\put(0,1){\circle*{0.4}}\put(0.2,1){\line(1,0){3}}
\put(3.4,1){\circle{0.4}}\put(3.6,1){\line(1,0){3}}\put(6.8,1){\circle{0.4}}\put(7,1){\line(1,0){1}}
\put(9,1){\circle*{0.2}} \put(10,1){\circle*{0.2}} \put(11,1){\circle*{0.2}} \put(12,1){\line(1,0){1}}
\put(13.2,1){\circle{0.4}}\put(13.4,1.1){\line(1,0){3}}\put(13.4,0.9){\line(1,0){3}}\put(14.6,1){\line(2,1){1}}
\put(14.6,1){\line(2,-1){1}}\put(16.6,1){\circle{0.4}} \put(16.6,1){\circle{0.6}} \put(0,0){1} \put(3.4,0){2}
\put(6.8,0){3} \put(13.2,0){$n-1$} \put(16.6,0){$n$}
\end{picture}\\
\hline $D_n$&
\begin{picture}(20,2.5)(-0.2,0)
\put(0,1){\circle*{0.4}} \put(0,1){\circle{0.6}} \put(0.2,1){\line(1,0){3}}
\put(3.4,1){\circle{0.4}}\put(3.6,1){\line(1,0){3}}\put(6.8,1){\circle{0.4}}\put(7,1){\line(1,0){1}}
\put(9,1){\circle*{0.2}} \put(10,1){\circle*{0.2}} \put(11,1){\circle*{0.2}} \put(12,1){\line(1,0){1}}
\put(13.2,1){\circle{0.4}}\put(13.4,1){\line(4,1){3}}\put(13.4,1){\line(4,-1){3}}\put(16.6,1.8){\circle*{0.4}}
\put(16.6,0.3){\circle*{0.4}} \put(16.6,1.8){\circle{0.6}} \put(16.6,0.3){\circle{0.6}} \put(0,0){1}
\put(3.4,0){2} \put(6.8,0){3} \put(12.5,0){$n-2$} \put(17,0){$n-1$} \put(17,1.5){$n$}
\end{picture}\\
\hline $E_6$&
\begin{picture}(20,2.5)(-0.2,0)
\put(0,0.7){\circle*{0.4}} \put(0,0.7){\circle{0.6}}
\put(0.2,0.7){\line(1,0){3}}\put(3.4,0.7){\circle{0.4}}\put(3.6,0.7){\line(1,0){3}}
\put(6.8,0.7){\circle{0.4}}\put(6.8,0.9){\line(0,1){1}}\put(6.8,2.1){\circle{0.4}}\put(7,0.7){\line(1,0){3}}
\put(10.2,0.7){\circle{0.4}}\put(10.4,0.7){\line(1,0){3}} \put(13.6,0.7){\circle*{0.4}}
\put(13.6,0.7){\circle{0.6}} \put(0.2,-0.1){1} \put(3.6,-0.1){3} \put(7,-0.1){4} \put(7,1.9){2}
\put(10.4,-0.1){5} \put(14,-0.1){6}
\end{picture}\\
\hline $E_7$&
\begin{picture}(20,2.5)(-0.2,0)
\put(0,0.7){\circle{0.4}} \put(0.2,0.7){\line(1,0){3}}\put(3.4,0.7){\circle{0.4}}\put(3.6,0.7){\line(1,0){3}}
\put(6.8,0.7){\circle{0.4}}\put(6.8,0.9){\line(0,1){1}}\put(6.8,2.1){\circle{0.4}}\put(7,0.7){\line(1,0){3}}
\put(10.2,0.7){\circle{0.4}}\put(10.4,0.7){\line(1,0){3}} \put(13.6,0.7){\circle{0.4}}
\put(13.8,0.7){\line(1,0){3}} \put(17,0.7){\circle*{0.4}} \put(17,0.7){\circle{0.6}} \put(0.2,-0.1){1}
\put(3.6,-0.1){3} \put(7,-0.1){4} \put(7,1.9){2} \put(10.4,-0.1){5} \put(13.8,-0.1){6} \put(17.5,-0.1){7}
\end{picture}\\
\hline
\end{tabular}\\
\begin{center}
\begin{picture}(25,2)(0,0)
\put(5,1){\circle*{0.4}} \put(6,0.8){minuscule} \put(12,1){\circle{0.4}} \put(12,1){\circle{0.6}}
\put(13,0.8){cominuscule} \put(21,0.8){both} \put(20,1){\circle*{0.4}} \put(20,1){\circle{0.6}}
\end{picture}
\end{center}

There are no minuscule nor cominuscule fundamental weight in type $E_8$, $F_4$, $G_2$.\\

\textbf{Assumption.} For the rest of the paper, the parabolic subgroup $P$ is assumed to be cominuscule.\\

\textbf{Acknowledgements.} I would like to thank Christian~Ohn for helpful discussions, and Takeshi~Ikeda for pointing out several references in the literature. 
I am also grateful to the referee for his valuable remarks, and especially for pointing out a gap in the proof of Proposition~\ref{intersection} and for providing a way 
to fill it.

\section{Local coordinates}\label{coordinates}

The notations are as in the Introduction. Moreover, $R_P$ denotes the root system associated with $P$:
\[
R^+ \setminus R_P^+ = \{ \beta \in R^+ \ |\ U_{\beta} \subset R_u(P) \},
\]
where $R_u(P)$ is the unipotent radical of $P$, and $U_{\beta}$ is the root subgroup associated with $\beta$.\\

Let $m \in X_w^v$. Then $m$ lies in a Schubert cell $C_{\tau}$ for some $\tau \in W^P$. Let
\[
U_{\tau}^-=\prod_{\beta \in \tau(R^+ \setminus R_P^+)}U_{-\beta}
\]
and $\mathcal{O}_{\tau}=U_{\tau}^-.e_{\tau}$, where $e_{\tau}=\tau P$. We identify $U_{-\beta}$ with $\C$ via an
isomorphism $\theta_{-\beta}: \C \to U_{-\beta}$ satisfying
\[
t\theta_{-\beta}(x)t^{-1}=\theta_{-\beta}\left(\frac{1}{\beta(t)}x\right)
\]
for all $ t \in T$ and all $x \in \C$. Let $N$ be the cardinality of $R^+ \setminus R_P^+$. We identify $\Ot$
with the affine space $\A^N$ via the isomorphism

\begin{equation}\label{isomorphism}
\begin{array}{rcl}
\A^N &\longrightarrow & \mathcal{O}_{\tau}\\
(x_{-\beta})_{\beta \in \tau(R^+ \setminus R_P^+)} & \mapsto & \displaystyle\prod_{\beta \in \tau(R^+ \setminus
R_P^+)} \theta_{-\beta}(x_{-\beta}).e_{\tau}.
\end{array}
\end{equation}
(In particular, $N$ is the dimension of $G/P$.)

\begin{lem} Let $\beta \in R$, and $\tau \in W^P$. Then $U_{\beta}$ fixes $e_{\tau}$ if and only if $-\beta \notin
\tau(R^+ \setminus R_P^+)$.
\end{lem}

\pr Let $\beta \in R$, and $\tau \in W^P$. Then
\[
\begin{split}
U_{\beta}.e_{\tau}=e_{\tau} & \iff \tau^{-1}U_{\beta}\ \tau P=P\\
                           & \iff U_{\tau^{-1}\beta} \subset P\\
                           & \iff \tau^{-1}\beta \in R^+ \ \textrm{or}\ -\tau^{-1}\beta \in R_P^+\\
                           & \iff -\beta \notin \tau(R^+) \ \textrm{or}\ -\beta \in \tau(R_P^+)\\
                           & \iff -\beta \notin \tau(R^+ \setminus R_P^+).~\square
\end{split}
\]

\begin{lem}\label{Ct} The Schubert cell $C_{\tau}$ is the affine subspace of $\mathcal{O}_{\tau}$ defined by the vanishing
of the coordinates $x_{-\beta}$ with $\beta \in R^+$.
\end{lem}

\pr Since $B$ is the semi-direct product of $T$ and the unipotent subgroup $U$, we have $C_{\tau}=U.e_{\tau}$.
Moreover, for any ordering of positive roots $\{\beta_1,\dots,\beta_p\}$,
\[
U=\prod_{i=1}^p U_{\beta_i}.
\]
We choose an ordering such that the positive roots $\beta$ with $-\beta \notin \tau(R^+ \setminus R_P^+)$ appear
at the end. Then, by the preceding lemma, we have:
\[
C_{\tau}=\prod_{\substack{\beta \in \tau(R^+ \setminus R_P^+)\\{}\beta <0}} U_{-\beta}.e_{\tau} \subset
\mathcal{O}_{\tau}.~\square
\]

The following lemma will be useful for the next section.

\begin{lem}\label{commute} For all $\beta, \gamma \in \tau(R^+ \setminus R_P^+)$ and for all $x, y \in \C$, the elements $\theta_{\beta}(x)$ and
$\theta_{\gamma}(y)$ commute.
\end{lem}

\pr We use the following expansion for the commutator (\textit{cf.}\  \cite{S}, proposition 8.2.3):
\[
\theta_{\beta}(x)\theta_{\gamma}(y)\theta_{\beta}(x)^{-1}\theta_{\gamma}(y)^{-1}=
\prod_{\substack{i\beta+j\gamma \in R\\ i,j>0}} \theta_{i\beta+j\gamma}(c_{\beta,\gamma,i,j}\ x^i\ y^j),
\]
where $c_{\beta,\gamma,i,j}$ are some constants in $\C$. Since the commutator must lie in $U_{\tau}^-$, it
suffices to prove that the roots of the form $i\beta+j\gamma$ do not lie in $\tau(R^+ \setminus R_P^+)$. Now, $P$ 
is the parabolic subgroup associated with the simple root $\a$. Since $\alpha$ is cominuscule, a positive root $\delta$ 
lies in $R^+ \setminus R_P^+$ if and only if $\alpha$ occurs with coefficient~1 in the expression of $\delta$. 
Clearly, $\alpha$ occurs with a coefficient $i+j$ in $\tau^{-1}(i\beta+j\gamma)$.~$\square$\\

\begin{Rem}~\label{iso}
 Identifying $\Ot$ with $U_\t^-$, it follows from Lemma~\ref{commute} that the isomorphism of algebraic varieties (\ref{isomorphism}) $\A^N \to \Ot$ is also an isomorphism of unipotent groups.
\end{Rem}

\begin{Ex}\label{SLn}
Let $G=SL_n(\C)$. It is a group of type $A_{n-1}$. The torus $T$ is the group of diagonal matrices of
determinant 1, and the Borel subgroup $B$ is the group of upper triangular matrices of determinant 1. The roots
are denoted $\alpha_{i,j}$, where
\[
\alpha_{i,j}:  T \to \C^*: \left(\begin{array}{cccc}
                      t_1 &  &  & \\
                      & t_2 &  &\\
                      & & \ddots & \\
                      &&&t_n\\
                      \end{array}\right)  \mapsto  \frac{t_i}{t_j}.
\]
The positive roots are the $\alpha_{i,j}$ with $i <j$, and the simple roots are the $\alpha_i=\alpha_{i,i+1}$
($i=1,\dots,n-1$). Let $\omega=\omega_d$ be the fundamental weight associated with the simple root $\alpha_d$.
The corresponding parabolic subgroup $P$ is
\[
P=\left\{\left(\begin{array}{c|c}
          * & *\\
          \hline
          0_{(n-d)\times d} & *\\
          \end{array}\right)\right\}.
\]
The group $G$ acts transitively on the Grassmannian $G_{d,n}$ of $d$-spaces in $\C^n$, and $P$ is the isotropy
subgroup of the vector space generated by $e_1,\dots,e_d$, where $(e_1,\dots,e_n)$ is the canonical basis of
$\C^n$. The Weyl group $W$ of this root system is $S_n$, and $W_P$ is isomorphic to $S_d \times S_{n-d}$, so
\[
W^P=I_{d,n}=\{\ui=i_1\dots i_d\ |\ 1 \leq i_1 <i_2< \dots < i_d \leq n\}.
\]
The Lie algebra $\mathfrak{g}$ of $G$ is the space of traceless matrices. Let $\mathfrak{t}$ be the Lie algebra
of the torus $T$. We have the weight decomposition of $\mathfrak{g}$:
\[
\mathfrak{g}=\mathfrak{t} \oplus \bigoplus_{i\neq j} \C E_{i,j}
\]
where $E_{i,j}$ is the elementary matrix with a 1 on the row $i$ and column $j$, and zero elsewhere. Thus, the
root subgroups are given by
\[
U_{\alpha_{i,j}}=\{I_n+x E_{i,j} \ |\ x \in \C \}
\]
and the isomorphism $\theta_{\alpha_{i,j}}$ is just $x \mapsto \exp(x E_{ij})$. Moreover,
\[
R^+ \setminus R_P^+=\{\alpha_{i,j} \ |\ i \leq d < j\},
\]
so in this case, Lemma~\ref{commute} becomes an elementary matrix computation.
\end{Ex}
~\\

Returning to the general case, we denote by $(m_{-\beta}|\beta \in \tau(R^+ \setminus R_P^+))$ the coordinates
of $m$, that is,
\[
m=\prod_{\beta \in \tau(R^+ \setminus R_P^+)} \theta_{-\beta}(m_{-\beta}).e_{\tau}.
\]

\begin{Nots}
We set:
\[
Y_w=X_w \cap \mathcal{O}_{\tau}, \quad Y^v=X^v \cap \mathcal{O}_{\tau},\quad Y_w^v=X_w^v \cap
\mathcal{O}_{\tau}.
\]
These sets are affine varieties, \textit{i.e.} Zariski-closed in $\mathcal{O}_{\tau}=\A^N$.
\end{Nots}

We now investigate if these affine varieties are cones over $m$.

\begin{prop}\label{Ywcone}
The varieties $Y_w$, $Y^v$ and $Y_w^v$ are cones over $e_{\tau}$.
\end{prop}

\pr Let $\o^\vee: \C^* \to T$ be the fundamental coweight associated to $P$. Since $\o^\vee$ is minuscule, the pairing $\langle \o^\vee , \g \rangle$ 
is equal to $1$ if $\g \in R^+ \setminus R_P^+$ (and to 0 if $\g \in R_P^+$). Now multiplication in $\A^N$ by a scalar $\xi$ is then given by conjugation 
in $U_\t^-$ by $\t(\o^\vee)(\xi)^{-1} \in T$: indeed, for $\b=\t(\g)$ with $\g \in R^+ \setminus R_P^+$, and for $z \in \C$, we
have
\begin{equation}\label{conjug}
 \t(\o^\vee)(\xi)^{-1} \theta_{-\beta} (z) \t(\o^\vee)(\xi)=\theta_{-\beta}(\xi^{\langle \t(\o^\vee),\b \rangle}z)=\theta_{-\beta}(\xi^{\langle \o^\vee,\g \rangle}z)=\theta_{-\beta}(\xi z).
\end{equation}

Let $x \in Y_w$ (resp. $x \in Y^v$), and $(x_{-\beta})$ be its coordinates. Then the point that has coordinates $(\xi x_{-\beta})$ is $t.x$, 
where $t=\t(\o^\vee)(\xi) \in T$. Therefore, this point lies in $X_w \cap \Ot$ (resp. in $X^v \cap \Ot$),
since $X_w$ (resp. $X^v$) is $T$--stable. It follows that $Y_w$, $Y^v$, and therefore $Y_w^v$ are cones over $e_\t$.~$\square$

\begin{prop}\label{Ywconem}
The variety $Y_w$ is a cone over $m$.
\end{prop}

\pr Consider the translation that maps $e_{\tau}$ to $m$. It is given in coordinates by $(x_{-\beta}) \mapsto
(x_{-\beta}+m_{-\beta})$. But if $x$ has coordinates $(x_{-\beta})$, then, by Remark~\ref{iso} the point of coordinates
$(x_{-\beta}+m_{-\beta})$ corresponds to $b.x$, where 
$b=\prod_{\beta} \theta_{-\beta}(m_{-\beta})$. Since $m_{-\beta}=0$ for all $\beta
>0$, we have $b \in B$ according to Lemma~\ref{Ct}. Now $b$ leaves $Y_w$ invariant and maps $e_{\tau}$ to $m$.~$\square$\\

However, the opposite Schubert variety $Y^v$ need not be a cone over $m$.

\begin{Ex}
We take the same notations as in Example~\ref{SLn}. In particular, using the identification $W^P=I_{d,n}$, we denote a Schubert variety in $G_{d,n}$ by $X_{i_1 \dots i_d}$,
and similarly for opposite Schubert and Richardson varieties. 
In the Grassmannian $G_{3,7}$, consider the Richardson variety $X_{356}^{125}$. The coordinates on the
open set $\mathcal{O}_{256}$ are parametrized by the set 
$\{12,15,16,32,35,36,42,45,46,72,75,76\}$ where $ij$ stands for the root $\alpha_{i,j}$. More precisely, we
have:

\[
\begin{array}{rcc}
 \A^{12} & \longrightarrow & \mathcal{O}_{256}\\
 (x_{12},x_{15},\dots,x_{76}) & \mapsto & \left[\begin{array}{ccc}
                                                  x_{12} & x_{15} & x_{16}\\
                                                  1 & 0 & 0\\
                                                  x_{32} & x_{35} & x_{36}\\
                                                  x_{42} & x_{45} & x_{46}\\
                                                  0 & 1 & 0\\
                                                  0 & 0 & 1\\
                                                  x_{72} & x_{75} & x_{76}
                                                  \end{array}\right].
\end{array}
\]
Here, a matrix between brackets actually stands for the 3-space in $\C^7$ generated by its columns. The
equations of $X_{356}$ are:
\[
\left\{\begin{array}{l}
        x_{72}=x_{75}=x_{76}=0\\
        x_{42}=0
        \end{array}
        \right.
\]
The equations of $X^{125}$ are:
\[
\left\{\begin{array}{l}
       x_{15}x_{36}-x_{35}x_{16}=0\\
       x_{15}x_{46}-x_{45}x_{16}=0\\
       x_{35}x_{46}-x_{45}x_{36}=0
       \end{array}\right.
\]
Let
\[
m=\left[\begin{array}{ccc}
1 & 0 & 1\\
1 & 0 & 0\\
0 & 0 & -1\\
0 & 0 & 0\\
0 & 1 & 0\\
0 & 0 & 1\\
0 & 0 & 0
\end{array}\right] \in X_{356}^{125}.
\]

We set:
\[
\left\{\begin{array}{ll}
y_{16}=x_{16}-1&\\
y_{36}=x_{36}+1&\\
y_{ij}=x_{ij} & \text{if\ } ij \notin \{16,36\}
\end{array}\right.
\]
The equations in these new coordinates are:
\[
\left\{\begin{array}{l}
        y_{72}=y_{75}=y_{76}=0\\
        y_{42}=0
        \end{array}
        \right.
\]
for $X_{356}$ and
\[
\left\{\begin{array}{l}
       y_{15}(y_{36}-1)-y_{35}(y_{16}+1)=0\\
       y_{15}y_{46}-y_{45}(y_{16}+1)=0\\
       y_{35}y_{46}-y_{45}(y_{36}-1)=0
       \end{array}\right.
\]
for $X^{125}$. While the equations for $X_{356}$ remain homogeneous, those for $X^{125}$ do not.
\end{Ex}
~\\

If $Y^v$ is indeed a cone over $m$, then we have the following result. The proof is taken from \cite{KL},
Remark~7.6.6.

\begin{prop}
Assume $Y^v$ is a cone over $m$. Let $\mu_w$ (resp.\ $\mu^v$, $\mu_w^v$) be the multiplicity of $m$ on $X_w$
(resp.\ $X^v$, $X_w^v$). Then
\begin{equation}
\mu_w^v=\mu_w\ \mu^v.
\end{equation}
\end{prop}

\pr In this case, $Y_w^v=Y_w \cap Y^v$ is a cone (over $m$) as well, so we may consider the projective varieties
$\PP(Y_w)$, $\PP(Y^v)$ and $\PP(Y_w^v)$, consisting of lines through $m$. Then $\mu_w$ (resp.\ $\mu^v$,
$\mu_w^v$) is just the degree of $\PP(Y_w)$ (resp.\ $\PP(Y^v)$, $\PP(Y_w^v)$). We conclude with B\'ezout's theorem
since $\PP(Y_w)$ and $\PP(Y^v)$
intersect transversely (\textit{cf.}\ \cite{R}, Corollary~1.5).~$\square$\\

\begin{Hyp}\label{notcone}
For the rest of the paper, we assume that $Y^v$ is not a cone over $m$.
\end{Hyp}

It is not clear however whether $Y_w^v$ is a cone or not. This problem will be solved in Section~\ref{proof}.

\section{Central projection and proof of Theorem~\ref{theformula}}\label{projection}

We shall compute the multiplicity of a point $m\in Y_w^v$ by relating it to degrees of projections, which requires us to work in a projective setting. More precisely, embed $\A^N$ into $\PP^N$ via

\[
\begin{array}{rrcl}
\iota: & \A^N & \hookrightarrow & \PP^N=\{[\xi:x_{-\beta}]\}\\
        & (x_{-\beta}) & \mapsto & [1:x_{-\beta}]
\end{array}
\]
and consider the projective closures
\[
Z_w=\overline{\iota(Y_w)},\qquad Z^v=\overline{\iota(Y^v)},\qquad Z_w^v=\overline{\iota(Y_w^v)}.
\]

We also identify $\PP^{N-1}$ with the hyperplane at infinity $\xi=0$ and consider the central projection
$p_m:\PP^N\to\PP^{N-1}$, sending any point $x \neq m$ to the intersection of the line $(mx)$ with $\PP^{N-1}$.
If $X\subset\PP^N$ is any projective variety and $m\in X$, then we have the following formula (\textit{cf.}\ \cite{M},
Theorem~5.11):
\begin{equation}\label{formulamult}
\!\!\deg X-\mult_mX=\begin{cases}
     \deg(p_m)_{|X}\,\deg(p_mX) & \textrm{\hspace*{-2mm}if $X$ is not a cone over $m$,}\\
     0 & \textrm{\hspace*{-2mm}if $X$ is a cone over $m$,}
     \end{cases}
\end{equation}
where $\deg X$ is the degree of $X$, $\deg(p_m)_{|X}$ is the degree of the rational map $p_m$ restricted to $X$, and $p_mX$ denotes the Zariski closure of $p_m(X \setminus \{m\})$.

\begin{prop}\label{theproposition}~
\begin{enumerate}
\renewcommand{\theenumi}{\alph{enumi}}
\item $\deg Z_w^v=\deg Z_w\deg Z^v$.
\item $Z^v_w$ is not a cone over $m$.
\item $\deg (p_m)_{|Z_w^v}=\deg(p_m)_{|Z^v}$.
\item $\deg (p_m Z_w^v)=\deg Z_w\deg (p_mZ^v)$.
\end{enumerate}
\end{prop}
We defer the proof to Section~\ref{proof}.\\

\emph{Proof of Theorem~\ref{theformula}.}
Using (\ref{formulamult}) and Proposition~\ref{theproposition}, we obtain
\[\begin{split}
\mu_{w}^{v}&=\deg Z_{w}^{v}-\deg (p_m)_{|Z_w^v}\,\deg (p_m Z_{w}^{v})\\
             &=\deg Z_{w}\deg Z^{v}-\deg (p_m)_{|Z^v}\,\deg Z_{w}\deg (p_m Z^{v})\\
             &=\deg Z_{w}\bigl[\deg Z^v-\deg (p_m)_{|Z^v}\,\deg (p_m Z^v)\bigr]\\
             &=\mu_{w}\ \mu^{v}.\ \square
\end{split}\]

\begin{Rem}\label{smooth}
In particular, this result enables us to find the singular locus of $X_w^v$ in terms of those of $X_w$ and $X^v$:
the point $m$ is smooth on $X_w^v$ if and only if $\mu_w^v=1$ if and only if $\mu_{w}=\mu^{v}=1$, that is, if
and only if $m$ is smooth on both $X_w$ and $X^v$. Note that this may also be seen more directly, using the fact
that $X_w$ and $X^v$ intersect properly and transversely at any point at which $\mu_{w}=\mu^{v}=1$
(\textit{cf.}\ \cite{R} Corollary~1.5, or \cite{BC} Corollary~2.9).
\end{Rem}

\section{$\C$-action on $G/P$}\label{Caction}

In this section, we introduce the main tool that will permit us to prove Proposition~\ref{theproposition} in the
next section. Let $e_\tau,m\in\Ot$ be as before: we shall construct an action of (the additive group) $\C$ on $G/P$ 
for which $e_\t$ and $m$ are in the same orbit.

Consider first the map 
\[
\begin{array}{lrcl}
 \phi^*: & \C^* & \to & B\\
         & \xi & \mapsto & \phi_\xi=\t(\o^\vee)(\xi)^{-1}b\t(\o^\vee)(\xi),
\end{array}
\]
where $b \in B \cap \U_\t^-$ is the element defined in the proof of Proposition~\ref{Ywconem}. The computation~(\ref{conjug}) shows that this 
map extends to a group homomorphism $\phi: \C \to B$. The natural $B$-action on $G/P$ thus induces a $\C$-action:
\[
 \begin{array}{lrcl}
  \Phi: & \C \times G/P & \to & G/P\\
        & (\xi,x) & \mapsto & \phi_{-\xi}.x
 \end{array}
\]

Moreover, $\Ot$ is invariant under this action (again by (\ref{conjug})). Actually, $\C$ acts on $\Ot=\A^N$ by translations: indeed, we get 
the following commutative diagram

\begin{equation}\label{diagramme}
\begin{array}{cc}
\xymatrix{\C \times \A^N \ar@{-->}[d]_{f}  \ar[r]^{\Phi} & \A^N \ar@{-->}[d]^{p_{e_{\tau}}}\\
            \PP^N \ar@{-->}[r]_{p_m} &  \PP^{N-1} } & \xymatrix{(\xi,x_{-\beta}) \ar@{|->}[d]_{f} \ar@{|->}[r]^{\Phi \qquad} & (x_{-\beta}-\xi m_{-\beta}) \ar@{|->}[d]^{p_{e_{\tau}}}\\
                                                     [\xi:x_{-\beta}] \ar@{|->}[r]_{p_m \qquad} & [0:x_{-\beta}-\xi m_{-\beta}]}
\end{array}
\end{equation}

Let us now restrict to $Y_w^v$: since it is a cone over $e_\tau$, a point $[\xi:x]$ lies in $Z_w^v$ if and only if $x \in Y_w^v$. It
follows that $f(\C \times Y_w^v)=Z_w^v$. Thus, the commutative diagram (\ref{diagramme}) restricts to
\begin{equation}\label{diagrestreint}
\begin{array}{c}
\xymatrix{\C \times Y_w^v \setminus \{(\xi,\xi m_{-\beta})|\ \xi \in \C\} \ar@{->>}[d]_{f} \ar[r]^{\qquad \Phi} &
\Phi(\C \times
Y_w^v)\setminus \{e_{\tau}\} \ar[d]^{p_{e_{\tau}}}\\
Z_w^v \setminus \{m\} \ar[r]_{p_m} & \PP^{N-1}}
\end{array}
\end{equation}

\begin{Rem}\label{fibre}
 Since (\ref{diagramme}) is a fibre product diagram, any fibre $\Phi^{-1}(\lambda y)$ (for $\lambda \neq 0$ and $[y] \in \PP^{N-1}$)
is mapped isomorphically via $f$ to the fibre $p_m^{-1}([y])$. Since we have the equalities $f(\C \times Y_w)=Z_w$, $f(\C \times Y^v)=Z^v$, $f(\C \times Y_w^v)=Z_w^v$ and 
$\C \times Y_w=f^{-1}(Z_w)$, $\C \times Y^v=f^{-1}(Z^v)$, $\C \times Y_w^v=f^{-1}(Z_w)$, the fibres of $\Phi_{|\C \times Y_w}$, $\Phi_{|\C \times Y^v}$, 
$\Phi_{|\C \times Y_w^v}$ over a point $\lambda y$ are isomorphic to the fibres of $p_{m|Z_w}$, $p_{m|Z^v}$, $p_{m|Z_w^v}$ over the point $[y]$.
\end{Rem}

In the next section, this remark will allow us to relate the degree of $p_m$ in diagram~(\ref{diagrestreint}) to that of $\Phi$.

\section{Proof of Proposition~\ref{theproposition}}\label{proof}

{\bf Proof of (a).} Since $Y_w$, $Y^v$, and $Y_w^v$ are (affine) cones over $e_{\tau}$, it is clear that $Z_w^v=Z_w \cap Z^v$. In addition, this intersection is proper and generically transverse (\cite{R}, Corollary~1.5), hence
$\deg Z_w^v=\deg Z_w\deg Z^v$ by B\'ezout's theorem.

\begin{Nots}
We denote by $F_w^v$ the closure in $\A^N$ of $\Phi(\C \times Y_w^v)$, and by $d_w^v$ the degree of $p_m: Z_{w}^{v} \setminus \{m\} \to p_m Z_w^v$ whenever it makes
sense (\textit{i.e.} when $Z_w^v$ is not a cone). We define $F_w, F^v, d^v$ in a similar way.
\end{Nots}

\begin{prop}\label{dvw}
When defined, the degree $d_w^v$ is equal to the degree of $\Phi: \C \times Y_{w}^{v} \to F_w^v$.
\end{prop}

\pr This follows from Remark~\ref{fibre}.~$\square$ 

\begin{lem}\label{fibredim}
The following properties are equivalent:
\begin{itemize}
 \item $Z_w^v$ is a cone over $m$,
 \item $F_w^v=Y_w^v$,
 \item every fibre of $\Phi: \C \times Y_w^v \to F_w^v$ has dimension~1.
\end{itemize}
In particular, they are true for $v=e$, hence $F_w=Y_w=\Phi(\C \times Y_w)$.
\end{lem}

\pr By Remark~\ref{fibre}, we see that the dimension of a generic fibre of $\Phi$ equals the dimension of
a generic fibre of $p_m$. Now $Z_w^v$ is a cone over $m$ if and only if every fibre of $p_m$ has dimension~1, if
and only if $\dim F_w^v =\dim Y_w^v $. But $Y_w^v=\Phi(0 \times Y_w^v) \subset F_w^v$ and the varieties $Y_w^v$
and $F_w^v$ are irreducible, so $Z_w^v$ is a cone over $m$ if and only if $F_w^v=Y_w^v$.~$\square$
~\\

{\bf Proof of (b) and (c).} By Proposition~\ref{dvw}, it suffices to compare the degree $d^v$ of $\Phi^v: \C \times Y^v \to F^v$ with
the degree $d_w^v$ of $\Phi_w^v: \C \times Y_w^v \to F_w^v$. First, the fibre of a point $x \in G/P$ for $\Phi$
is
\[
\Phi^{-1}(x)=\{(\xi,\Phi(-\xi,x))\ |\ \xi \in \C\}.
\]
In particular, a point lies in $\im(\Phi^v)$ (resp.\ in $\im(\Phi_w^v)$) if and only if its $\C$-orbit meets
$Y^v$ (resp.\ $Y_w^v$). There exists an open set $\Omega^v$ of $F^v$ in which the fibre of every point $y$ consists
of $d^v$ points. Then $d^v$ is just the number of points in the $\C$-orbit
of $y$ that belong to $Y^v$. Now set $y=(y_{-\beta})_{\beta \in \tau(R^+ \setminus R_P^+)}$ and let
\[
c=\prod_{\substack{\beta \in \tau(R^+ \setminus R_P^+)\\ \beta <0}}\theta_{-\beta}(y_{-\beta})
\qquad
c^-=\prod_{\substack{\beta \in \tau(R^+ \setminus R_P^+)\\ \beta >0}}\theta_{-\beta}(-y_{-\beta}),
\]
so we have $c.e_{\tau}=c^-.y=:x$. Since $c \in B$, $x \in C_{\tau} \subset Y_w$. Now $c^-$ commutes with
$\phi_{\xi}$ for all $\xi \in \C$, hence every point in $c^-(\Omega^v)$ has a $\C$-orbit which
meets $Y^v$ in exactly $d^v$ points. In particular, $F_w^v \neq Y_w^v$, since otherwise every fibre of $\Phi_w^v$
would have dimension~1 (by Lemma~\ref{fibredim}), which is not the case for the fibre of $x$. This already shows
(b), so it makes sense to talk about the degree $d_w^v$ of $\Phi_w^v$. Thus, let $\Omega_w^v$ be an open set of
$F_w^v$ such that for every point $z$ in $\Omega_w^v$, the fibre of $z$ consists of $d_w^v$ points. Since $x \in
c^-(\Omega^v)$, $c^-(\Omega^v) \cap F_w^v$ and $\Omega_w^v$ are non-empty open sets of the irreducible variety
$F_w^v$, so they must meet. Taking $z$ in this intersection, we see that $d_w^v=d^v$, which shows
(c).~$\square$\\

\begin{prop}\label{interproper}
The intersection $F_w \cap F^v$ is proper and transverse on an open set of $F_w^v$.
\end{prop}

\pr The transversality of the intersection $F_w \cap F^v$ on a generic point in $F_w^v$ follows from the transversality of the
intersection of a direct Schubert variety and an opposite Schubert variety. More precisely, let $(F_w)_{sm}$ be the open set of smooth points
of $F_w$. Taking a point smooth on $Y_w^v$ shows that $\Omega_w=(F_w)_{sm}\cap F_w^v$ is a non-empty open set of $F_w^v$. Let $(F^v)_{sm}$ be the open set
of smooth points of $F^v$. Again, $\Omega^v=(F^v)_{sm} \cap F_w^v \neq \emptyset$. Indeed, take a smooth point $x$ of $F^v$ belonging to
$\Phi(\C \times Y^v)$. We have seen in  the previous proof that from $x$ we can construct an isomorphism $c^-$ of $F^v$
mapping $x$ to a point of $F_w^v$, which thus remains smooth on $F^v$. The two non-empty open subsets $\Omega_w$ and $\Omega^v$ of the irreducible
variety $F_w^v$ have a non-empty intersection $\Omega_w^v$. Now $O_w^v=\Phi^{-1}(\Omega_w^v) \cap (\PP^1 \times Y_w^v)_{sm} \neq \emptyset$
since $\PP^1 \times Y_w^v$ is irreducible. We claim that $\Phi: O_w^v \to \Omega_w^v$ is dominant. Indeed, we must show that every open subset $U$
of $\Omega_w^v$ meets $\Phi(O_w^v)$. Since $U$ is open in $F_w^v$, $U \cap \Phi(\C \times Y_w^v) \neq \emptyset$. So it makes sense to talk
about $\Phi^{-1}(U)$, which is an open set of $\C \times Y_w^v$. Thus, $\Phi^{-1}(U) \cap O_w^v \neq \emptyset$, which
implies $U \cap \Phi(O_w^v) \neq \emptyset$. Since $\Phi: O_w^v \to \Omega_w^v$ is dominant, we know that $\Phi(O_w^v)$ contains a non-empty open
set $\Omega$ of $\Omega_w^v$. Let us summarize the properties of $\Omega$: it is a non-empty open subset of $F_w^v$, whose every point $y$ is smooth
in both $F_w$ and $F^v$, and $y=\Phi(p)$ with $p$ smooth in $\C \times Y_w^v$, so $p$ is smooth in both $\C \times Y_w$ and $\C \times Y^v$.

Let $y=\Phi(p) \in \Omega$ be such a point. We view the map $\Phi:\C \times \A^N \to \A^N: (\xi,x) \mapsto
\phi_{-\xi}.x$ as a map $\Phi: \C^{N+1} \to \C^{N}$. It is linear and surjective. Thus,
\[
\begin{split}
\C^N\supset T_y(F_{w})+T_y(F^{v})&\supset d\Phi_p(T_p(\C \times Y_{w}))+d\Phi_p(T_p(\C \times Y^{v}))\\
                         &\supset d\Phi_p(\C \oplus (T_pY_{w}+T_pY^{v}))\\
                         &\supset d\Phi_p(\C \oplus \C^N)\\
                         &\supset \C^N.
\end{split}
\]

This transversality result proves that the intersection is proper: indeed, on one hand, $\dim(F_w \cap F^v) \ge \dim(F_w)+\dim(F^v)-N$,
but on the other hand,
\[
 \begin{split}
   \dim(F_w \cap F^v) & \le \dim(T_y(F_w \cap F^v)) \le \dim(T_yF_w \cap T_yF^v)\\
                      & \le \dim(T_yF_w)+\dim(T_yF^v)-\dim(T_yF_w+T_yF^v)\\
                      & \le \dim(F_w)+\dim(F^v)-N.\ \square
 \end{split}
\]
\\

\begin{prop}\label{intersection}
We have the equality $F_w^v=F_w \cap F^v$. In particular, the intersection $F_w \cap F^v$ is generically transverse.
\end{prop}

This result will be proved in the next section.\\

{\bf Proof of (d).} Since $y=\Phi(\xi,x)$ implies $zy=\Phi(z\xi,zx)$ for all $z\in\C$, $\Phi(\C\times Y_w^v)$ is a cone over $e_\tau$, and so is its closure $F_w^v$.
But by the commutative diagram (\ref{diagrestreint}),
\[
p_{e_\tau}(F_w^v \setminus \{e_\tau\}) \subset
\overline{p_{e_\tau}(\Phi(\C \times Y_w^v) \setminus
\{e_\tau\})}=p_m Z_w^v.
\]
Comparing dimensions, we see that $p_{e_\tau}F_w^v=p_mZ_w^v$, \emph{i.e.}\ $p_mZ_w^v$ is the projective variety at infinity of the cone $F_w^v$.
In particular, $\deg(p_mZ_w^v)=\deg(F_w^v)$, and similarly $\deg(p_mZ_w)=\deg(F_w)$ and $\deg(p_mZ^v)=\deg(F^v)$. Equality (d) now follows from Proposition~\ref{interproper}
and B\'ezout's theorem, noting that $\deg(p_mZ_w)=\deg(Z_w)$.~$\square$\\

\section{Proof of Proposition~\ref{intersection}} \label{secinter}

Since $\Phi(\C \times Y_w^v) \subset \Phi(\C \times Y_w) \cap \Phi(\C \times Y^v)$, we obtain $F_w^v \subset F_w \cap F^v$.
Moreover, the first inclusion is an equality: indeed, if $z=\Phi(\xi,x)\in Y_w$ with $\xi \in \C, x \in Y^v$, then $x=\Phi(-\xi,z)\in Y_w$ 
since $\Phi(\C \times Y_w)=Y_w$, so $z=\Phi(\xi,x) \in \Phi(\C \times Y_w^v)$.

However, the inclusion $F_w \cap F^v \subset F_w^v$ requires some work. Let $\mathcal{U}=\{(\xi,x,\Phi(\xi,x))|\xi \in \C, x \in
G/P\}$ and $\Gamma$ be its closure in $\PP^1 \times G/P \times G/P$ (so $\Gamma$ is the graph of $\Phi$ viewed as a rational map). We have a commutative diagram:
\[\begin{array}{cc}
\xymatrix{\Gamma \ar[d]_{\pi_1 \times \pi_2}\ar[dr]^{\pi_3} & \\
\PP^1 \times G/P \ar@{-->}[r]_{\Phi} & G/P} & \xymatrix{(\xi,x,y) \ar@{|->}[d]_{\pi_1 \times \pi_2}\ar@{|->}[dr]^{\pi_3} & \\
(\xi,x) \ar@{|->}[r]_{\Phi} & \Phi(\xi,x)} \end{array}\]

The morphism $\pi_1 \times \pi_2: \Gamma \to \PP^1 \times G/P$ is surjective, and restricts to an isomorphism
between $\mathcal{U}$ and $\C \times G/P$. In particular, $\Gamma$ is an irreducible projective variety of
dimension $N+1$. 

Likewise, let $\U_w=\{(\xi,x,\Phi(\xi,x))|\xi \in \C, x \in X_w\}$ and $\Gamma_w$ be its closure, and similarly for $\U^v$, $\U_w^v$, $\Gamma^v$, $\Gamma_w^v$. 
Then $\pi_3(\Gamma_w)=\pi_3(\overline{\U_w})=\overline{\pi_3(\U_w)}$ in $G/P$, so $\pi_3(\Gamma_w)\cap \Ot$ is the closure of $\pi_3(\U_w)\cap \Ot=\Phi(\C \times Y_w)$ in $\Ot$.
Proceeding similarly with $\Gamma^v$ and $\Gamma_w^v$, we obtain 
\[
\pi_3(\Gamma_w) \cap \Ot=F_w,\quad \pi_3(\Gamma^v) \cap \Ot=F^v,\quad \pi_3(\Gamma_w^v) \cap \Ot=F_w^v.
\]

We now need to study the $\pi_3$-fibre of a point in $F_w$. Actually, if $y$ is in $Y_w$, then its
fibre lies entirely in $\Gamma_w$. Indeed, $\Ut$ naturally acts on $G/P$ and on $\Gamma$ via
$g.(\xi,x,y)=(\xi,g.x,g.y)$ (since $\Ut$ is Abelian), and the morphism $\pi_3$ is $\Ut$-equivariant. It follows that
whenever two points in $G/P$ belong to the same $\Ut$-orbit, their fibres are isomorphic. Now since $\pi_3: \Gamma \to
G/P$ is dominant, there is an open set in $G/P$ in which every point has a fibre of pure
dimension~1. Since $\Ot$ is open in $G/P$, it meets this open set, and since $\Ot$ is a $\Ut$-orbit in $G/P$,
$y$ itself has a fibre of pure dimension~1. 

Now fix an irreducible component $C$ of $\pi_3^{-1}(y)$. Then
\[
\bigl(\pi_1 \times \pi_2 (C)\bigr)\cap(\C\times G/P)\subset\Phi^{-1}(y).
\]
If $C\cap\mathcal{U}\ne\emptyset$, then the left hand side of this inclusion is non-empty and of dimension~1.
Since $\Phi^{-1}(y)$ is isomorphic to the $\C$-orbit of $y$, it is itself irreducible of dimension (at most)~1,
hence the inclusion becomes an equality. Taking closures, we then obtain $C=\overline{\{(\xi,x,y) | (\xi,x) \in\Phi^{-1}(y)\}}$;
in particular, $C$ is the unique irreducible component of $\pi_3^{-1}(y)$ that intersects $\mathcal{U}$.
Note also that $C\subset \Gamma_w$.

Now let $C'$ be an irreducible component of $\pi_3^{-1}(y)$ different from $C$, so that $C' \subset \{\infty\} \times G/P \times \{y\}$.
Let $\Gamma_\infty \subset \Gamma$ be the subvariety $\pi_1^{-1}(\infty)$. We have a $\Ut$-equivariant morphism
$\pi: \Gamma_\infty \to G/P: (\infty,x,y) \mapsto y$, so $C'$ is an irreducible subvariety of the fibre $\pi^{-1}(y)$.
Since $\Gamma_\infty \subsetneq \Gamma$, its dimension is at most $N$. Because of the equivariance of $\pi$, we
see that $\Ot$ is in the image of $\pi$, so $\pi$ is surjective. Decomposing $\Gamma_\infty$ into irreducible
components $\Gamma_{\infty}=C_1 \cup \dots \cup C_r$, we obtain $G/P=\pi(C_1) \cup \dots \cup
\pi(C_r)$, so that for some $i$, $\pi: C_i \to G/P$ is onto. Renumbering the $C_i$, we may assume that for some $t \ge 1$, 
$C_1,\dots,C_t$ are mapped surjectively to $G/P$, and $C_{t+1},\dots,C_r$ are not. For $i \leq t$, there is an open 
set $U_i$ of $G/P$ such that each element on $U_i$ has a finite fibre in
$C_i$. For $i>t$, let $U_i$ be the open set $G/P \setminus \pi(C_i)$. Taking the intersection $U=\bigcap_{i=1}^n
U_i$, we obtain a non-empty open set of $G/P$ satisfying the following property: for each $z \in U$, the fibre
of $z$ in $\Gamma_\infty$ consists of a finite number of points. Again, $U$ meets the open orbit $\Ot$, so this
property is true for every point in $\Ot$, in particular for $y$. So $C'$ is included in the finite fibre $\pi^{-1}(y)$:
a contradiction. Therefore, $C'$ cannot exist, \emph{i.e.}\ $\pi_3^{-1}(y)=C\subset\Gamma_w$ is irreducible, and not contained 
in $\{\infty\} \times G/P \times G/P$.

Assume now that $F_w^v \neq F_w \cap F^v$. By Proposition~\ref{interproper}, $F_w^v$ and $F_w \cap F^v$ have the same dimension, thus $F_w \cap F^v$
is not irreducible. Let $F$ be an irreducible component of the intersection $F_w \cap F^v$ different from $F_w^v$. Let $y \in F$, and assume that 
$y \notin F_w^v$. Then $y \notin \pi_3(\mathcal{U}^v)$, so $\pi_3^{-1}\{y\} \subset \Gamma^v \setminus \mathcal{U}^v \subset \{\infty\} \times G/P \times G/P$.
But $y \in F_w$, and we have seen that in this case $\pi_3^{-1}(y)$ is never contained in $\{\infty\} \times G/P \times G/P$. This gives a contradiction.~$\square$

\appendix
\section*{Appendix.\ Singularities of Schubert varieties in $SO(2n+1)/P_1$}\label{SO}

In this Appendix, we shall determine the singular locus of Schubert varieties in $G/P$, where $G$ is of type
$B_n$ and $P$ is cominuscule. So let $V=\C^{2n+1}$ together with a non-degenerate symmetric bilinear form
$(.,.)$ given in the canonical basis $(e_1,\dots,e_{2n+1})$ by the anti-diagonal matrix $E$ with 1's all along
the anti-diagonal. The expression of the
quadratic form $Q$ associated with $(.,.)$ is
\[
Q(x_1,\dots,x_{2n+1})=x_{n+1}^2+2\sum_{i=1}^n x_ix_{2n+2-i}.
\]
Let $G=SO(V)$, $B \subset G$ the subgroup of upper triangular matrices, and $T \subset G$ the subgroup
of diagonal matrices. Then $B$ is a Borel subgroup of $G$ and $T$ is a maximal torus of $G$. The
group $G$ acts naturally on $V$, and $e_1$ is a highest weight vector, of weight $\omega_1$ (the unique
cominuscule weight of $G$), so that $G/P_1$ gets identified with the $G$-orbit of the line generated by $e_1$:
\[
G/P_1=\{[x_1:\dots:x_{2n+1}]\ |\ Q(x_1,\dots,x_{2n+1})=0\ \}
\]

In this setting, the Schubert varieties are given by
\[
X_i=\{[x_1:\dots:x_i:0:\dots:0]\ |\ Q(x_1,\dots,x_i,0,\dots,0)=0\},
\]
with $1\le i\le 2n+1$, but $i\ne n+1$. Indeed, let $x=[x_1:\dots:x_{i-1}:1:0:\dots:0]$
with $Q(x)=0$, and let us prove that $x \in C_i$, that is, there exists $b \in B$ such that $x=b.e_i$. A straightforward calculation shows that we may take the columns $b_1,\dots,b_{2n+1}$ of $b$ as follows:
\begin{itemize}
\item Case 1: $i<n+1$.\\[2mm]
\hspace*{1cm}$b_j=\left\{\begin{array}{ll}
              e_j & \text{if}\  j \neq i\ \text{and}\ j \leq 2n+2-i\\
              x & \text{if}\ j=i\\
              e_j-x_{2n+2-j}e_{2n+2-i} & \text{otherwise}
              \end{array}
              \right.$\\[3mm]
\item Case 2: $i>n+1$.\\[2mm]
\hspace*{1cm}$b_j=\left\{\begin{array}{ll}
              e_j & \text{if}\  j \leq 2n+2-i \\
              x & \text{if}\ j=i\\
              e_j-x_{2n+2-j}e_{2n+2-i} & \text{otherwise}
              \end{array}
              \right.$\\[3mm]
\end{itemize}

The Jacobian criterion easily shows that $\Sing X_i$ is equal to $X_{2n+1-i}$ if $i>n+1$, and empty if $i<n+1$. Moreover, since $X_i$ is defined by a single quadratic equation, the multiplicity of a singular point must be equal to~2. Hence there are two cases for the multiplicity $\mu_i(x)$ of a point $x=[x_1:\dots:x_i:0:\dots:0]$ on $X_i$:
\begin{itemize}
\item Case 1: $i<n+1$. Then $\mu_i(x)=1$.
\item Case 2: $i>n+1$. Then
\[
\mu_i(x)=\left\{\begin{array}{ll}
                       2 & \text{if $x_i=\dots=x_{2n+2-i}=0$,}\\
                       1 & \text{otherwise.}
                       \end{array}\right.
\]
\end{itemize}

Of course, we have the same result for the opposite Schubert varieties
\[
X^j=\{[0:\dots:0:x_j:\dots:x_{2n+1}]\ |\ Q(0,\dots,0,x_j,\dots,x_{2n+1})=0 \}
\]
There are again two cases for the multiplicity $\mu^j(x)$ of $x=[0:\dots:0:x_j:\dots:x_{2n+1}]$ on $X^j$:
\begin{itemize}
\item Case 1: $j<n+1$. Then
\[
\mu^j(x)=\left\{\begin{array}{ll}
                        2 & \text{if $x_j=\dots=x_{2n+2-j}=0$,}\\
                        1 & \text{otherwise.}
                        \end{array}\right.
\]
\item Case 2: $j>n+1$. Then $\mu^j(x)=1$.
\end{itemize}

Note that a Richardson variety $X_i^j$ ($j \leq i$) also is a quadric in a projective space, so the
multiplicity of a point $m \in X_i^j$ must be at most~2. But by Theorem~\ref{theformula}, if $m$ were singular
in both $X_i$ and in $X^j$, then its multiplicity would be~4. This means that $\Sing X_i\cap\Sing X^j=\emptyset$,
a fact that can also be verified directly: indeed, if this intersection is non-empty, then $2n+3-j \leq 2n+1-i$,
so $j \le i \le j-2$, a contradiction.\\

\end{document}